\documentclass[10pt]{article}

\newif\ifpdf
\ifx\pdfoutput\undefined
    \pdffalse       
\else
    \pdfoutput=1    
    \pdftrue
\fi

\usepackage{amsmath,amssymb,amsthm}

\ifpdf
    \usepackage[pdftex]{graphicx}
    \usepackage[backref, colorlinks=true, pdfstartview=FitH, linkcolor=blue,
        citecolor=blue, urlcolor=blue]{hyperref}
\else
    \usepackage{graphicx}
    \usepackage{hyperref}
\fi

\title{Sturmian Words \\and the \\Permutation that Orders Fractional Parts}
\author{Kevin O'Bryant \\ University of California, San Diego \\ kobryant@math.ucsd.edu \\ http://www.math.ucsd.edu/$\sim$kobryant}
\date{\today}

    \newcommand{\N}{{\mathbb N}}
    \newcommand{\Z}{{\mathbb Z}}
    \newcommand{\Q}{{\mathbb Q}}
    \newcommand{\R}{{\mathbb R}}
    \newcommand{\C}{{\mathbb C}}

    \newcommand{\floor}[1]{\left\lfloor #1 \right\rfloor}
    \newcommand{\ceiling}[1]{\big\lceil #1 \big\rceil}
    \newcommand{\fp}[1]{\{ #1 \}}
    \newcommand{\tf}[1]{[\![ #1 ]\!]}

    \DeclareMathOperator{\ord}{ord}
    
    \DeclareMathOperator{\sgn}{sgn}

    \newcommand{\vw}{\vec{w}^\sigma}
    \newcommand{\vwp}{\vec{w}^{\pi}}
    \newcommand{\vd}{\vec{\delta}}

    \newcommand{\al}{\alpha}
    \newcommand{\Ls}{{\cal L}_\sigma}
    \newcommand{\Ms}{{\cal M}_\sigma}
    \newcommand{\Ma}{{\cal M}_n(\alpha)}
    \newcommand{\Ps}{{\cal P}_\sigma}
    \newcommand{\pan}{\pi_{\alpha,n}}
    \newcommand{\B}{B_{\alpha}}
    \newcommand{\BB}{B^\prime}
    \newcommand{\QQ}{{\cal Q}}
    \newcommand{\II}{{\cal I}}
    \newcommand{\VV}{{\cal V}}

    \newtheorem{thm}{Theorem}[section]
    \newtheorem{lem}{Lemma}[subsection]
    \newtheorem{prop}[lem]{Proposition}
    \newtheorem{cor}[lem]{Corollary}

    \newenvironment{xmp}{\bigskip \noindent{\bf Example:}\;}{\hfill $\square$\bigskip}
    \newenvironment{cnj}{\bigskip \noindent{\bf Conjecture:}\;}{\bigskip}

    \newenvironment{proofof}[1]{\medskip\noindent{\em Proof of #1.}}{\qed\medskip}

\begin{document}
    \ifpdf
       \DeclareGraphicsExtensions{.pdf,.jpg,.mps,.png}
    \fi
\maketitle

\begin{abstract}
A Sturmian word is a map \mbox{$W:\N \to \{0,1\}$} for which the set of $\{0,1\}$-vectors
\linebreak\mbox{$F_n(W):=\{(W(i),W(i+1),\dots,W(i+n-1))^T:i\in\N\}$} has cardinality exactly
$n+1$ for each positive integer $n$. Our main result is that the volume of the simplex whose
$n+1$ vertices are the $n+1$ points in $F_n(W)$ does not depend on $W$. Our proof of this
motivates studying algebraic properties of the permutation $\pan$ (where $\al$ is any
irrational and $n$ is any positive integer) that orders the fractional parts
$\fp{\al},\fp{2\al},\dots,\fp{n\al}$, i.e.,
\mbox{$0<\fp{\pan(1)\al}<\fp{\pan(2)\al}<\dots<\fp{\pan(n)\al}<1$.} We give a formula for
the sign of $\pan$, and prove that for every irrational $\al$ there are infinitely many $n$
such that the order of $\pan$ (as an element of the symmetric group $S_n$) is less than $n$.
\end{abstract}

\section{Introduction}\label{sec.Introduction}

A {\em binary word} is a map from the nonnegative integers into $\{0,1\}$. The {\em factors
of $W$} are the column vectors $(W(i),W(i+1),\dots, W(i+n-1))^T$, where $i\ge0$ and $n\ge1$.
In particular, the set of factors of length $n$ of a binary word $W$ is defined by
 $$F_n(W) := \left\{ \left(W(i), W(i+1),\dots, W(i+n-1)\right)^T : i\ge0\right\}.$$
Obviously, $|F_n(W)| \le 2^n$ for any binary word $W$. It is known
\cite[Theorem~1.3.13]{2002.Lothaire} that if $|F_n(W)|<n+1$ for any $n$, then $W$ is
eventually periodic. If $|F_n(W)|=n+1$ for every $n$---the most simple non-periodic
case---then $W$ is called a {\em Sturmian word}. Sturmian words arise in many fields,
including computer graphics, game theory, signal analysis, diophantine approximation,
automata, and quasi-crystallography. The new book of Lothaire~\cite{2002.Lothaire} provides
an excellent introduction to combinatorics on words; the second chapter is devoted to
Sturmian words.

Throughout this paper, $W$ is always a Sturmian word, $n$ is always a positive integer, and
$\al$ is always an irrational between 0 and 1. A typical example of a Sturmian word is given
by $c_{\al}(i):=\floor{(i+2)\al}-\floor{(i+1)\al}$, the so-called {\em characteristic word
with slope $\al$}. By routine manipulation, one finds that $c_{\al}(i)=1$ if and only if
$i+1\in \{ \floor{ k\al}\colon k\in\Z^+\}$. The integer sequences
$(\floor{k\al+\beta})_{k=1}^\infty$ are called {\em Beatty sequences}. The study of Beatty
sequences is intimately related to the study of Sturmian words, and the interested reader
can locate most of the literature through the bibliographies of \cite{Stolarsky1976},
\cite{Brown1993}, and \cite{Tijdeman2000a}.

In this paper, we consider the $n+1$ factors in $F_n(W)$ to be the vertices of a simplex in
$\R^n$. Our main result is
\begin{thm}\label{thm.Volume}
If $W$ is a Sturmian word, then the volume of the simplex $F_n(W)$ is $\frac 1{n!}$.
\end{thm}
The remarkable aspect of Theorem~\ref{thm.Volume} is that the volume of the simplex $F_n(W)$ is
independent of $W$. The key to the proof of Theorem~\ref{thm.Volume} is to study $F_n(W)$ for
all Sturmian words $W$ simultaneously. The primary tool is the representation theory of finite
groups.

Sturmian words are examples of one-dimensional quasicrystals, at least with respect to some
of the `working definitions' currently in use. In contrast to the study of crystals, group
theory has not been found very useful in the study of quasicrystals. According to M.
Senechal~\cite{SenechalBook}, ``The one-dimensional case suggests that symmetry may be a
relatively unimportant feature of aperiodic crystals.'' Thus, the prominent role of
symmetric groups in the proof of Theorem~\ref{thm.Volume} comes as a surprise.

The proof of Theorem~\ref{thm.Volume} reveals a deep connection between the simplex
$F_n(c_\al)$ and algebraic properties of the permutation $\pan$ of $1,2,\dots,n$ that orders
the fractional parts $\fp{\al}, \fp{2\al}, \dots, \fp{n\al}$, i.e.,
 $$0<\fp{\pan(1)\al} < \fp{\pan(2)\al} < \dots < \fp{\pan(n)\al} < 1.$$

The definition of $\pan$ has a combinatorial flavor, and accordingly some attention has been
given to its combinatorial qualities. Using the geometric theory of continued fractions,
S{\'o}s~\cite{Sos1957} gives a formula for $\pan$ in terms of $n$, $\pan(n)$, and $\pan(1)$ (see
Lemma \ref{lem.Sos}). Boyd \& Steele~\cite{1979.Boyd.Steele} reduce the problem of finding the
longest increasing subsequence in $\pan$ to a linear programming problem, which they then solve
explicitly. Schoi{\ss}engeier~\cite{Schoissengeier1984} used Dedekind eta sums to study $\pan$ and
give his formula for the star-discrepancy of $n\al$-sequences.

Here, motivated by the appearance of $\pan$ in our study of the simplex $F_n(W)$, we initiate
the study of algebraic properties of $\pan$. If $\sigma$ is an element of a group (with
identity element {\tt id}), we let $\ord(\sigma)$ be the least positive integer $t$ such that
$\sigma^t={\tt id}$, or $\infty$ if no such integer exists. We use this notation with
permutations, matrices, and congruence classes (the class will always be relatively prime to
the modulus). For any permutation $\sigma$, let $\sgn(\sigma)$ be the sign of $\sigma$, i.e.,
$\sgn(\sigma)=1$ if $\sigma$ is an even permutation and $\sgn(\sigma)=-1$ if $\sigma$ is an odd
permutation. Our main results concerning $\pan$ are stated in Theorems~\ref{thm.Order}
and~\ref{thm.Sign}.
\begin{thm}\label{thm.Order}
For every irrational $\al$, there are infinitely many positive integers $n$ such that
$\ord(\pan)<n$.
\end{thm}
\begin{thm}\label{thm.Sign}
For every irrational $\al$ and positive integer $n$,
 $$\sgn(\pi_{\al,2n})=\sgn(\pi_{\al,2n+1}) = \prod_{\ell=1}^n (-1)^{\floor{2\ell\al}}.$$
\end{thm}
In particular, although $\pan$ is ``quasi-random'' in the sense of~\cite{2002.Cooper}, it is
highly structured in an algebraic sense.

Sections~\ref{sec.Sturmian} and~\ref{sec.FracParts} are logically independent and may be read
in either order. In Section~\ref{sec.Sturmian}, we consider Sturmian words and the simplex
$F_n(W)$. Section~\ref{sec.FracParts} is devoted to proving Theorems~\ref{thm.Order} and
~\ref{thm.Sign}. Section~\ref{sec.Questions} is a list of questions raised by the results of
Sections~\ref{sec.Sturmian} and~\ref{sec.FracParts} that we have been unable to answer. A {\em
Mathematica} notebook containing code for generating the functions and examples in this paper
is available from the author.

\section{Sturmian Words}\label{sec.Sturmian}

\subsection{Introduction to Sturmian Words}\label{sec.Sturmian.sub.Introduction}

An excellent introduction to the theory of Sturmian words is given in
\cite[Chapter~2]{2002.Lothaire}. We restate the results needed in this paper in this
subsection.

If $\al\in(0,1)$ is irrational and $\beta$ is any real number, then the words $s_{\al,\beta}$
and $s_{\al,\beta}^\prime$ defined by
    \begin{align*}
        s_{\al,\beta} (i)          &:= \floor{(i+1)\al+\beta}-\floor{i\al+\beta} \\
        s_{\al,\beta}^\prime (i)    &:= \ceiling{(i+1)\al+\beta}-\ceiling{i\al+\beta}
    \end{align*}
are Sturmian, and every Sturmian word arises in this way \cite[Theorem 2.3.13]{2002.Lothaire}.
The irrational number $\al$ is called the slope of the word, and the word $c_\al :=
s_{\al,\al}$ is called the characteristic word of slope $\al$. It is easily shown
\cite[Proposition 2.1.18]{2002.Lothaire} that $F_n(W)$ depends only on the slope of $W$, and so
it is consistent to write $F_n(\al)$ in place of $F_n(W)$. In fact, we shall use the equation
$F_n(\al)=F_n(s_{\al,\beta})$ for all $\beta$. It is often easier to think in terms of `where
the 1s are'; elementary manipulation reveals that
    $$c_\al(i) = \begin{cases} 1 & i+1 \in \{ \floor{\tfrac{k}{\al}} : k\ge1\} \\ 0 &
    \text{otherwise.} \end{cases}
    $$

The $n+1$ elements of $F_n(\alpha)$ are $n$-dimensional vectors, naturally defining a simplex
in $\R^n$. Whenever a family of simplices arises, there are several questions that must be
asked. Can the simplex $F_n(\alpha)$ be degenerate? If $F_n(\alpha)$ is not degenerate, can one
express its volume as a function of $n$ and $\alpha$? Under what conditions on $\alpha, \beta,
n$ is $F_n(\alpha)\cong F_n(\beta)$?

The first and second questions are answered by Theorem~\ref{thm.Volume}, which we prove in
Subsection~\ref{sec.Sturmian.sub.Proof} below. Computer calculation suggests a simple answer to
the third question, which we state as a conjecture in Section~\ref{sec.Questions}.

\begin{xmp}
The characteristic word with slope $e^{-1}\approx 0.368$ begins
    $$(c_{e^{-1}}(0), c_{e^{-1}}(1), c_{e^{-1}}(2), \dots)
        = (0, 1, 0, 0, 1, 0, 0, 1, 0, 1, 0, 0, 1, 0, 0, 1, 0, 0, 1, 0, 1, \dots).$$
Note that
    $$c_{e^{-1}}(i) = \begin{cases} 1 & i+1 \in \{\floor{n e}\colon n\ge 1\} = \{2,5,8,10,\dots\} \\
    0 & \text{otherwise.} \end{cases}$$
The set of factors of $c_{e^{-1}}$ of length 6, arranged in anti-lexicographic order, is
    $$F_6(c_{e^{-1}}) = F_6(e^{-1}) = \left\{
        \begin{pmatrix}
            1\\0\\1\\0\\0\\1
        \end{pmatrix},
        \begin{pmatrix}
             1\\0\\0\\1\\0\\1
        \end{pmatrix},
        \begin{pmatrix}
             1\\0\\0\\1\\0\\0
        \end{pmatrix},
        \begin{pmatrix}
             0\\1\\0\\1\\0\\0
        \end{pmatrix},
        \begin{pmatrix}
             0\\1\\0\\0\\1\\0
        \end{pmatrix},
        \begin{pmatrix}
             0\\0\\1\\0\\1\\0
        \end{pmatrix},
        \begin{pmatrix}
             0\\0\\1\\0\\0\\1
        \end{pmatrix} \right\}.$$
\end{xmp}

\subsection{Definitions}\label{sec.Sturmian.sub.Definitions}

To analyze a simplex, one first orders the vertices (we order them anti-lexicographically).
Then, one translates the simplex so that one vertex is at the origin (we move the last factor
to $\vec{0}$). Finally, one writes the coordinates of the other vertices as the columns of a
matrix. If this matrix is non-singular, then the simplex is not degenerate. In fact, the volume
of the simplex is the absolute value of the determinant divided by $n!$. We are thus led to
define the matrix $\Ma$, whose $j$-th column is $\vec{v}_{j}-\vec{v}_{n+1}$, where
$F_n(\alpha)=\{\vec{v}_1,\vec{v}_2,\dots,\vec{v}_{n+1}\}$ ordered anti-lexicographically.

\begin{xmp}
$${\cal M}_6(e^{-1}) =
\begin{pmatrix}
 1 & 1  & 1  & 0  & 0  & 0  \\
 0 & 0  & 0  & 1  & 1  & 0  \\
 0 & -1 & -1 & -1 & -1 & 0  \\
 0 & 1  & 1  & 1  & 0  & 0  \\
 0 & 0  & 0  & 0  & 1  & 1  \\
 0 & 0  & -1 & -1 & -1 & -1 \\
\end{pmatrix}.$$
\end{xmp}

When a list of vectors is enclosed by parentheses, it denotes a matrix whose first column is
the first vector, second column the second vector, and so on. For example,
 $$\Ma := \left( \vec{v}_1-\vec{v}_{n+1}, \vec{v}_2-\vec{v}_{n+1}, \dots
 \vec{v}_n-\vec{v}_{n+1}\right).$$
We also define $\VV_k$ to be the $n\times n$ matrix all of whose entries are 0, save the
$k$-th column, which is $\vec{e}_{k-1}-2\vec{e}_k+\vec{e}_{k+1}$.

We shall make frequent use of Knuth's notation:
 $$\tf{Q}=\begin{cases} 1 & \text{$Q$ is true;} \\ 0 & \text{$Q$ is false.}\end{cases}$$

We denote the symmetric group on the symbols $1,2,\dots,n$ by $S_n$. We use several notations
for permutations interchangeably. We use standard cycle notation when convenient, and
frequently use {\em one-line notation} for permutations, i.e.,
    $$\sigma = [\sigma(1),\dots,\sigma(n)].$$
Thus, if a list of distinct numbers is surrounded by parentheses then it is a permutation in
cycle notation, and if the numbers $1, 2, \dots, n$ are in any order and surrounded by brackets
then it is a permutation in one-line notation. We multiply permutations from right to left.
Also, set $\Ps=(p_{ij})$, with $p_{ij}=\tf{j=\sigma(i)}$. This is the familiar representation
of $S_n$ as permutation matrices.

One permutation we have already defined is $\pan$. For notational convenience we set
$\pan(0):=0$ and $\pan(n+1):=n+1$. Also, set $P_j:=\fp{j\al}$ for $0\le j \le n$, and set
$P_{n+1}:=1$. Thus
 $$0=P_{\pan(0)} < P_{\pan(1)} < P_{\pan(2)} < \dots < P_{\pan(n)} < P_{\pan(n+1)}=1.$$

We write $\vec{e}_i$ ($1\leq i \leq n$) be the $n$-dimensional column vector with every
component 0 {\em except} the $i$-th component, which is 1. We set $\vec{e}_{n+1}=\vec{0}$,
the $n$-dimensional 0 vector. We denote the identity matrix as
$\II:=(\vec{e}_1,\vec{e}_2,\dots,\vec{e}_n)$. Let $\vec{\delta}_i:=\vec{e}_{i+1}-\vec{e}_i$
($1\leq i \leq n$). In particular, $\vec{\delta}_n=-\vec{e}_n$.

We will also use the notation $h(\vec{v})$ for the Hamming weight of the $\{0,1\}$-vector
$\vec{v}$, i.e., the number of 1's.

Set $$D(\sigma):=\{1\} \cup \{k \colon \sigma^{-1}(k-1)>\sigma^{-1}(k)\}.$$ In other words,
$D(\sigma)$ consists of those $k$ for which $k-1$ does not occur before $k$ in
$[\sigma(1),\sigma(2),\dots,\sigma(n)]$. For example, $D([1,3,5,4,2,6])=\{1,3,5\}$.

Set $$\vw_1 :=\sum_{i\in D(\sigma)} \vec{e}_i$$ and for $1\leq j \leq n$, set
$$\vw_{j+1}:=\vw_j+\vd_{\sigma(j)}.$$

We now define two matrices: the $n\times (n+1)$ matrix
 $$\Ls:=(\vw_1,\vw_2,\dots,\vw_n,\vw_{n+1}),$$
and the square $n\times n$ matrix
 $$\Ms:=(\vw_1-\vw_{n+1},\vw_2-\vw_{n+1},\dots,\vw_n-\vw_{n+1}).$$
Proposition~\ref{prop.GoodDefinitions} below shows that $\Ma = {\cal M}_{\pan}$, justifying
our definitions.

\begin{xmp}
Set $n=5$ and $\sigma=[5,2,3,1,4]=(1,5,4)(2)(3)$. We find that $D(\sigma)=\{1,2,5\}$, and so
$\vw_1=\vec{e}_{1}+\vec{e}_2+\vec{e}_5$. By definition $\vw_2=\vw_1+\vd_{\sigma(1)} =
\vw_1+\vec{e}_6-\vec{e}_5 = \vec{e}_1+\vec{e}_2$, and so on. Thus
 $$\Ls=\begin{pmatrix}
 1 & 1 & 1 & 1 & 0 & 0 \\
 1 & 1 & 0 & 0 & 1 & 1 \\
 0 & 0 & 1 & 0 & 0 & 0 \\
 0 & 0 & 0 & 1 & 1 & 0 \\
 1 & 0 & 0 & 0 & 0 & 1
    \end{pmatrix}
\quad \text{and} \quad
 \Ms=\begin{pmatrix}
 1 & 1  & 1  & 1  & 0  \\
 0 & 0  & -1 & -1 & 0  \\
 0 & 0  & 1  & 0  & 0  \\
 0 & 0  & 0  & 1  & 1  \\
 0 & -1 & -1 & -1 & -1
    \end{pmatrix}.$$
Note that the first column of $\Ms$ is $\vec{e}_1$; that this is always the case is proven in
Lemma~\ref{lem.2.4.2}. Further, the second column of $\Ms$ is $\vec{e}_1+\vd_{\sigma(1)}$, the
third is $\vec{e}_1+\vd_{\sigma(1)}+\vd_{\sigma(2)}$, and so forth. This pattern holds in
general and is proved in Lemma~\ref{lem.2.4.3} below. It is not immediate from the definitions
that $\Ls$ is always a $\{0,1\}$-matrix or that $\Ms$ is a $\{-1,0,1\}$-matrix; we prove this
in Lemma~\ref{lem.2.4.4}.

In Lemma~\ref{lem.map.is.1-1} we prove that if $\sigma\not=\tau$ then $\Ms\not={\cal
M}_{\tau}$. The proof relies on reconstructing $\sigma$ and $\Ls$ from $\Ms$. This
reconstruction proceeds as follows. The `$-1$' entries of $\Ms$ are in the second and fifth
rows; this gives $\vw_{6}=\vec{e}_2+\vec{e}_5$, which is the last column of $\Ls$. In fact,
the $j$-th column of $\Ls$ is the $j$-th column of $\Ms$ plus $\vec{e}_2+\vec{e}_5$. Once we
know the columns of $\Ls:=\left(\vw_1,\vw_2,\dots,\vw_6\right)$, we can use the definition
of $\vw_{j+1}$ to find $\sigma(j)$. For example,
$\vd_{\sigma(4)}=\vw_5-\vw_4=\vec{e}_2-\vec{e}_1=\vd_{1}$, and so $\sigma(4)=1$.

Lemma~\ref{lem.M_transposition} generalizes the observation that
 $${\cal M}_{[1,2,4,3,5]}={\cal M}_{(4,3)}=
   \begin{pmatrix}
 1 & 0 & 0 & 0 & 0 \\
 0 & 1 & 0 & 0 & 0 \\
 0 & 0 & 1 & 1 & 0 \\
 0 & 0 & 0 & -1 & 0 \\
 0 & 0 & 0 & 1 & 1
    \end{pmatrix} = \II+\VV_4.$$

With $\phi=\frac{\sqrt 5 -1}{2}$, we compute that $\pi_{\phi,5}=[5,2,4,1,3]=(1,5,3,4)(2)$,
and one may directly verify that ${\cal M}_{\pi_{\phi,5}}={\cal M}_{\phi}(5)$. This is no
accident, by Proposition~\ref{prop.GoodDefinitions} below ${\cal M}_{\alpha}(n)={\cal
M}_{\pan}$ for all $\alpha$ and $n$. The equation $${\cal M}_{(4,3)} {\cal M}_{\pi_{\phi,5}}
= {\cal M}_{(4,3)(1,5,3,4)(2)}={\cal M}_{(1,5,4)(2)(3)},$$ which is the same as
  $$
   \begin{pmatrix}
 1 & 0 & 0 & 0 & 0 \\
 0 & 1 & 0 & 0 & 0 \\
 0 & 0 & 1 & 1 & 0 \\
 0 & 0 & 0 & -1 & 0 \\
 0 & 0 & 0 & 1 & 1
    \end{pmatrix}
    \begin{pmatrix}
   1 & 1 & 1 & 1 & 0 \\
   0 & 0 & -1& -1& 0 \\
   0 & 0 & 1 & 1 & 1 \\
   0 & 0 & 0 & -1& -1 \\
   0 & -1& -1& 0 & 0
 \end{pmatrix}
 =
 \begin{pmatrix}
 1 & 1  & 1  & 1  & 0  \\
 0 & 0  & -1 & -1 & 0  \\
 0 & 0  & 1  & 0  & 0  \\
 0 & 0  & 0  & 1  & 1  \\
 0 & -1 & -1 & -1 & -1
    \end{pmatrix},$$
is an example of the isomorphism of Proposition~\ref{prop.representation}.
\end{xmp}

\subsection{The Matrices $\Ma$ and ${\cal M}_{\pan}$}\label{sec.Sturmian.sub.Matrices}

\begin{prop}\label{prop.GoodDefinitions}
$\Ma ={\cal M}_{\pan}$.
\end{prop}

\begin{proof}
For brevity, we write $\pi$ in place of $\pan$. First observe that
$\vwp_1,\vwp_2,\dots,\vwp_{n+1}$ are in anti-lexicographic order by definition, and so for
\mbox{$1\le i < j \le n+1$} we have $\vwp_i \not=\vwp_j$. We know from \cite[Proposition
2.1.18]{2002.Lothaire} that $F_n(\al)=F_n(s_{\al,\beta})$ for every $\beta$, and from
\cite[Theorem 2.1.13]{2002.Lothaire} that $|F_n(\al)|=n+1$. Thus, it suffices to show that
$\vwp_j \in F_n(s_{\al,\beta})$ for some $\beta$. In fact, we shall show that
    $$\left( s_{\al,\beta_j}(1), s_{\al,\beta_j}(2), \dots, s_{\al,\beta_j}(n)\right)
        =\vwp_j$$
with $\beta_j := -P_1-P_{\pi(j)}$.

Using the identities $\floor{x}=x-\fp{x}$ and $\fp{x-y}=\fp{x}-\fp{y}+\tf{\fp{x}<\fp{y}}$, we
have
 \begin{align}\label{eq.s.TF.Formula}
 s_{\al,\beta_j}(i)
    &=  \floor{(i+1)\al-P_1-P_{\pi(j)}}-\floor{i\al-P_1-P_{\pi(j)}} \notag\\
    &=  \al-P_i+P_{i-1}-\tf{P_i<P_{\pi(j)}}+\tf{P_{i-1}<P_{\pi(j)}} \notag\\
    &=  \tf{P_i<P_{i-1}}-\tf{P_i<P_{\pi(j)}}+\tf{P_{i-1}<P_{\pi(j)}}.
 \end{align}
The last equality follows from the knowledge that $s_{\al,\beta_j}(i)\in\Z$, and consequently
if $P_i<P_{i-1}$ then $\al-P_i+P_{i-1}>\al>0$ must in fact be 1, and if $P_i>P_{i-1}$ then
$\al-P_i+P_{i-1}<\al<1$ must in fact be 0.

We first consider $j=1$. We have $P_1>P_0$, $P_1 \ge P_{\pi(1)}$, and $P_0<P_{\pi(1)}$, whence
$s_{\al,\beta_1}(1)=1=\tf{1\in D(\pi)}$. For $2\le i \le n$, we have $P_i \ge P_{\pi(1)}$ and
$P_{i-1}\ge P_{\pi(1)}$, whence
 $$
 s_{\al,\beta_1}(i)
    =  \tf{P_i < P_{i-1}}
    =  \tf{\pi^{-1}(i) < \pi^{-1}(i-1)}
    =  \tf{i\in D(\pi)}.
 $$
Therefore, $\left( s_{\al,\beta_1}(1), s_{\al,\beta_1}(2), \dots, s_{\al,\beta_1}(n)\right) =
\vwp_1$.

Now suppose that $2\le j \le n+1$. Since $\vwp_j$ is defined by
$\vwp_{j}-\vwp_{j-1}=\vd_{\pi(j-1)}=\vec{e}_{\pi(j-1)+1}-\vec{e}_{\pi(j-1)}$, we need to
show that $s_{\al,\beta_j}(i)-s_{\al,\beta_{j-1}}(i)=\tf{i=\pi(j-1)+1} - \tf{i=\pi(j-1)}.$
By Eq.~(\ref{eq.s.TF.Formula}), we have
    \begin{align*}
    s_{\al,\beta_j}(i)-s_{\al,\beta_{j-1}}(i)
        &=  -\tf{P_i<P_{\pi(j)}}+\tf{P_{i-1}<P_{\pi(j)}}
                +\tf{P_i<P_{\pi(j-1)}}-\tf{P_{i-1}<P_{\pi(j-1)}} \\
        &=  \left(\tf{P_{i-1}<P_{\pi(j)}}-\tf{P_{i-1}<P_{\pi(j-1)}} \right)
                -\left(\tf{P_i<P_{\pi(j)}}-\tf{P_i<P_{\pi(j-1)}}\right) \\
        &=  \tf{i-1=\pi(j-1)} - \tf{i=\pi(j-1)}.
    \end{align*}
\end{proof}

We remark that in the same manner one may prove that if
 $$1-\fp{\pan(j)\al} \leq \fp{i\al+\beta}< 1-\fp{\pan(j-1)\al},$$
then
 $$\left( s_{\al,\beta}(i), s_{\al,\beta}(i+1), \dots, s_{\al,\beta}(i+n-1)\right) =
    \vec{w}_j^{\pan}.$$
From this it is easy to prove that $F_n(s_{\al,\beta})$ does not depend on $\beta$ and has
cardinality $n+1$, the two results of \cite{2002.Lothaire} that we used. There is another
fact that can be proved in this manner that we do not use explicitly but which may help the
reader develop an intuitive understanding of the simplices $F_n(W)$. If $r,r^\prime$ are
consecutive Farey fractions of order $n+1$ and $r<\al<\gamma<r^\prime$, then
$F_n(\al)=F_n(\beta)$.

\subsection{A Curious Representation}\label{sec.Sturmian.sub.Representation}

\begin{prop}\label{prop.representation}
The map $\sigma \mapsto {\cal M}_\sigma$ is an isomorphism.
\end{prop}

The proof we give of this is more a verification than an explanation. We remark that there are
several anti-isomorphisms involved in our choices. We have chosen to multiply permutations
right-to-left rather than left-to-right. We have chosen to consider the list $[a,b,c,\dots]$ as
the permutation taking 1 to $a$, 2 to $b$, 3 to $c$, etc., rather than the permutation taking
$a$ to 1, $b$ to 2, $c$ to 3, etc. Finally, we have chosen to define the vectors $\vw$ to be
columns rather than rows. If we were to change any two of these conventions, then we would
still get an isomorphism.

We begin with some simple observations about $\Ms$ before proving
Proposition~\ref{prop.representation}. We defined $\vw_{j+1}:=\vw_j+\vd_{\sigma(j)}$. An easy
inductive consequence of this definition is that $\vw_{j+1}=\vw_1+\sum_{i=1}^j \vd_{\sigma(i)}$
for $1\leq j \leq n$; we use this repeatedly and without further fanfare.

\begin{lem}\label{lem.2.4.2}
$\vw_1-\vw_{n+1}=\vec{e}_1$.
\end{lem}

\begin{proof}
Since $\vw_{n+1}=\vw_1+\sum_{i=1}^n \vd_{\sigma(i)}$, all we need to show is that $\sum_{i=1}^n
\vd_{\sigma(i)}=-\vec{e}_1$. As $\{1,2,\dots,n\}=\{\sigma(1),\sigma(2),\dots,\sigma(n)\}$, we
have
 $$\vw_{n+1}-\vw_1=\sum_{i=1}^n \vd_{\sigma(i)}=\sum_{i=1}^n \vd_i=\sum_{i=1}^n
(\vec{e}_{i+1}-\vec{e}_i) = \vec{e}_{n+1}-\vec{e}_1= -\vec{e}_1.\hfill \qed$$
\renewcommand{\qed}{}\end{proof}

\begin{lem}\label{lem.2.4.3}
The $j$-th column of $\Ms$ is $\vec{e}_1+\sum_{i=1}^{j-1} \vd_{\sigma(i)}$.
\end{lem}

\begin{proof}
The $j$-th column of $\Ms$ is defined as $\vw_j-\vw_{n+1}$. If $j=1$, then this lemma reduces
to Lemma~\ref{lem.2.4.2}. If $j>1$, then Lemma~\ref{lem.2.4.2} gives
 \begin{equation*}\vw_j-\vw_{n+1} =\left(\vw_1+\sum_{i=1}^{j-1} \vd_{\sigma(i)}\right)-
\left(\vw_1-\vec{e}_1\right) = \vec{e}_1+\sum_{i=1}^{j-1} \vd_{\sigma(i)}.\hfill \qed
 \end{equation*}
\renewcommand{\qed}{}\end{proof}

\begin{lem}\label{lem.2.4.4}
The entries of $\Ls$ are 0 and 1. The entries of $\Ms$ are -1, 0, and 1.
\end{lem}

\begin{proof}
It is easily seen from Lemma~\ref{lem.2.4.3} that $\Ms$ is a $\{-1,0,1\}$-matrix, and this also
follows from the more subtle observation that $\Ls$ is a $\{0,1\}$-matrix. To prove that $\Ls$
is a $\{0,1\}$-matrix is the same as showing that each $\vw_j$ ($1\leq j \leq n+1$) is a
$\{0,1\}$-vector. This is obvious for $\vw_1:= \sum_{i\in D(\sigma)} \vec{e}_i$. We have
$\vw_{j+1}=\vw_1+\sum_{i=1}^j \vd_{\sigma(i)}$; define $\vec{v}, c_i$ by $\vec{v}:=\sum_{i=1}^j
\vd_{\sigma(i)}=\sum_{i=1}^n c_i\vec{e}_i$.

The $i$-th component of $\vec{v}$ can only be affected by $\vd_{i-1}$ (which adds 1 to the
$i$-th component) and $\vd_i$ (which subtracts 1). It is thus clear that $c_i$ is (-1), (0),
(0), or (1) depending, respectively, on whether (only $i$), ($i-1$ and $i$), (neither $i-1$ nor
$i$), or (only $i-1$) is among $\{\sigma(1),\sigma(2),\dots,\sigma(j)\}$. To show that
$\vw_{j+1}=\vw_1+\vec{v}$ is a $\{0,1\}$-vector, we need to show two things. First, if $c_i=-1$
then $i\in D(\sigma)$ (and so the $i$-th component of $\vw_1$ is 1). Second, if $c_i=1$ then
$i\not\in D(\sigma)$ (and so the $i$-th component of $\vw_1$ is 0).

If $c_i=-1$, then $i$ is and $i-1$ is not among $\{\sigma(1),\dots,\sigma(j)\}$. This means
that $\sigma^{-1}(i-1)>j \geq \sigma^{-1}(i)$, and so by the definition of $D$, $i \in
D(\sigma)$. If, on the other hand, $c_i=1$, then $i-1$ is and $i$ is not among
$\{\sigma(1),\dots,\sigma(j)\}$. This means that $\sigma^{-1}(i-1)\leq j < \sigma^{-1}(i)$, and
so by the definition of $D$, $i \not\in D(\sigma)$.
\end{proof}

\begin{lem}\label{lem.map.is.1-1}
The map $\sigma\mapsto \Ms$ is 1-1.
\end{lem}

\begin{proof}
Given $\Ms$, we find $\sigma$. We first note that moving from $\Ls$ to $\sigma$ is easy since
$\vw_{j+1}-\vw_j=\vd_{\sigma(j)}$. Some effort is involved in finding $\Ls$ from $\Ms$. We need
to find $\vw_{n+1}$.

Every row of $\Ls$ contains at least one 0 (explanation below). Thus $\Ms$ will contain a $-1$
in exactly those rows in which $\vw_{n+1}$ contains a 1, and we are done. We have
$\vw_{j+1}=\vw_j+\vd_{\sigma(j)}=\vw_j+\vec{e}_{\sigma(j)+1}-\vec{e}_{\sigma(j)}$. The
$\sigma(j)$-th row of $\vec{e}_{\sigma(j)+1}$ is 0, and $\vw_j \,, \vec{e}_{\sigma(j)}$ are
$\{0,1\}$-vectors, so the $\sigma(j)$-th row of $\vw_{j+1}$ is either $0$ or $-1$. But by
Lemma~\ref{lem.2.4.4}, $\Ls$ is a $\{0,1\}$-matrix, whence the $\sigma(j)$-th row of the
$(j+1)$-st column of $\Ls$ is 0.
\end{proof}

Recall that $\VV_k$ is defined to be the matrix all of whose entries are 0, save the $k$-th
column, which is $\vec{e}_{k-1}-2\vec{e}_k+\vec{e}_{k+1}$.

\begin{lem}\label{lem.M_transposition}
${\cal M}_{(k,k-1)}=\II+\VV_k$, for $2\leq k \leq n$.
\end{lem}

\begin{proof}
Follow the definitions. With $\sigma=(k,k-1)$, we have $D(\sigma)=\{1,k\}$,
$\vec{w}^{(k,k-1)}_j=\vec{e}_j+\vec{e}_k$ for $n\geq j\not=k$,
$\vec{w}^{(k,k-1)}_k=\vec{e}_{k-1}+\vec{e}_{k+1}$ and $\vec{w}^{(k,k-1)}_{n+1}=\vec{e}_k$.
\end{proof}

We have laid the necessary groundwork, and turn now to proving
Proposition~\ref{prop.representation}.

\bigskip \begin{proof} We have already seen in Lemma~\ref{lem.map.is.1-1} that the map
$\sigma \mapsto \Ms$ is 1-1; all that remains is to show that this map respects multiplication,
i.e., for any $\sigma,\tau \in S_n$, ${\cal M}_\tau \Ms = {\cal M}_{\tau\sigma}$. Since we may
write $\tau$ as a product of transpositions of the form $(k,k-1)$ it is sufficient to show that
${\cal M}_{(k,k-1)}\Ms={\cal M}_{(k,k-1)\sigma}$ for every $k$ ($2\leq k\leq n$) and $\sigma
\in S_n$.

We need to split the work into two cases: $\sigma^{-1}(k-1)<\sigma^{-1}(k)$ and
$\sigma^{-1}(k-1)>\sigma^{-1}(k)$. In each case we first describe the rows of ${\cal
M}_{(k,k-1)}\Ms-\Ms$ using Lemma~\ref{lem.M_transposition}, and then compute the columns of
${\cal M}_{(k,k-1)\sigma}-\Ms$ from the definition. We will find that ${\cal
M}_{(k,k-1)}\Ms-\Ms={\cal M}_{(k,k-1)\sigma}-\Ms$ in each case, which concludes the proof.
Since the two cases are handled similarly, we present only the first case.

Suppose that $\sigma^{-1}(k-1)<\sigma^{-1}(k).$ By Lemma~\ref{lem.M_transposition}, ${\cal
M}_{(k,k-1)}\Ms-\Ms=\VV_k\Ms$. The matrix $\VV_k$ is zero except in the $(k-1,k)$, $(k,k)$, and
$(k+1,k)$ positions (note: we sometimes refer to positions which do not exist for $k=n$; the
reader may safely ignore this detail), where it has value 1, $-2$, 1, respectively. Thus
$\VV_k\Ms$ is zero except in the $(k-1)$-st and $(k+1)$-st rows (which are the same as the
$k$-th row of $\Ms$), and the $k$-th row (which is $-2$ times the $k$-th row of $\Ms$).

We now describe the $k$-th row of $\Ls$. By the hypothesis of this case, $k \not\in D(\sigma)$,
so the $k$-th row of $\vw_1$ is 0. Since $\vw_{j+1}=\vw_1+\sum_{i=1}^j \vd_{\sigma(i)}$, the
$k$-th row of $\vw_j$ is 0 for $1\leq j \leq \sigma^{-1}(k-1)$, is 1 for $\sigma^{-1}(k-1) < j
\leq\sigma^{-1}(k)$, and is 0 for $\sigma^{-1}(k) < j \leq n+1$. This gives the $k$-th row of
$\Ms$ as $\sigma^{-1}(k-1)$ `0's followed by $\sigma^{-1}(k)-\sigma^{-1}(k-1)$ `1's, followed
by $n-\sigma^{-1}(k)$ `0's.

We now compute ${\cal M}_{(k,k-1)\sigma}-\Ms$. The columns of ${\cal L}_{(k,k-1)\sigma}$ are
given by $\vec{w}^{(k,k-1)\sigma}_{j+1}= \vec{w}^{(k,k-1)\sigma}_1+\sum_{i=1}^j
\vd_{(k,k-1)\sigma(i)}$ (except the first, but the first column of ${\cal M}_\tau$ is
$\vec{e}_1$, independent of $\tau$). Now $(k,k-1)\sigma(i)=\sigma(i)$ for
$i\not\in\{\sigma^{-1}(k),\sigma^{-1}(k-1)\}$, so that $\sum_{i=1}^j
\vd_{(k,k-1)\sigma(i)}=\sum_{i=1}^j \vd_{\sigma(i)}$ for $j<\sigma^{-1}(k-1)$ and for $j\geq
\sigma^{-1}(k)$. For $\sigma^{-1}(k-1)\leq j < \sigma^{-1}(k)$,
 $$\sum_{i=1}^j \vd_{(k,k-1)\sigma(i)}
    =\left(\sum_{i=1}^j \vd_{\sigma(i)}\right)-\vd_{k-1}+\vd_{k}.$$
Thus the $(j+1)$-st column of ${\cal M}_{(k,k-1)\sigma}-\Ms$ is
 $$\left(\vec{w}^{(k,k-1)\sigma}_{j+1}-\vec{w}^{(k,k-1)\sigma}_{n+1}\right)-
 \left(\vw_{j+1}-\vw_{n+1}\right)
    = \sum_{i=1}^{j} \vd_{(k,k-1)\sigma(i)}-\sum_{i=1}^j \vd_{\sigma(i)},$$
which is $\vec{0}$ for $j<\sigma^{-1}(k-1)$ and for $j\geq \sigma^{-1}(k)$, and
$-\vd_{k-1}+\vd_k=\vec{e}_{k-1}-2\vec{e}_k+\vec{e}_{k+1}$ for $\sigma^{-1}(k-1)\leq j <
\sigma^{-1}(k)$. We have shown that ${\cal M}_{(k,k-1)}\Ms-\Ms={\cal M}_{(k,k-1)\sigma}-\Ms$ in
the case $\sigma^{-1}(k-1)<\sigma^{-1}(k)$.
\end{proof}

\subsection{Proof of Theorem \ref{thm.Volume}}\label{sec.Sturmian.sub.Proof}

\begin{proofof}{Theorem \ref{thm.Volume}}
The volume of the $n$-dimensional simplex whose vertices have coordinates $\vec{v}_1,
\vec{v}_2, \dots, \vec{v}_{n+1}$ is $\frac{1}{n!}$ times the absolute value of the determinant
of the matrix
    $$\left(
    \vec{v}_1-\vec{v}_{n+1},\vec{v}_2-\vec{v}_{n+1}, \dots, \vec{v}_{n}-\vec{v}_{n+1}
    \right).$$
In our case, this means that the volume of the simplex $F_n(\al)$ is $\tfrac{1}{n!}\left|
\det(\Ma) \right|$. We will show that $\det(\Ma)=\pm1$.

Now $\Ma={\cal M}_{\pan}$ by Proposition~\ref{prop.GoodDefinitions}, and for any integer $t$ we
have $\left({\cal M}_{\pan}\right)^t={\cal M}_{\pan^t}$ by
Proposition~\ref{prop.representation}. Since $\pan\in S_n$, a finite group, there is a positive
integer $t$ such that $\pan^t$ is the identity permutation (which we denote by {\tt id}). Thus
    \begin{equation*}
    \left(\det \Ma \right) ^t
    =   \left(\det {\cal M}_{\pan}\right)^t
    =   \det \left({\cal M}_{\pan}^t\right)
    =   \det \left({\cal M}_{\pan^t}\right)
    =   \det \left({\cal M}_{\tt id}\right)
    =   \det \II
    =   1.
    \end{equation*}
Consequently, $\det \Ma$ is a $t$-th root of unity, and since the entries of $\Ma$ are
integers, $\det \Ma = \pm1$.
\end{proofof}

\subsection{The Character of the Representation}\label{sec.Sturmian.sub.Character}

We review the needed facts and definitions from the representation theory of finite groups
(an excellent introduction is \cite{2001.Sagan}). A {\em representation} of a finite group
$G$ is a homomorphism ${\cal R}:G \to SL_m(\C)$ for some $m\geq1$. The representation is
said to be {\em faithful} if the homomorphism is in fact an isomorphism. Thus,
Proposition~\ref{prop.representation} implies that $\{\Ms \colon  \sigma \in S_n\}$ is a
faithful representation of $S_n$. The {\em character} of ${\cal R}$ is the map $g \mapsto
tr({\cal R}(g))$, with $tr$ being the trace function. We will use Corollary 1.9.4(5) of
\cite{2001.Sagan}, which states that if two representations ${\cal R}_1,{\cal R}_2$ of $G$
have the same character, then they are similar, i.e., there is a matrix $\QQ$ such that
$\forall g \in G \left(\QQ^{-1}{\cal R}_1(g)\QQ={\cal R}_2(g)\right)$. Any such matrix
${\cal Q}$ is called an {\em intertwining matrix} for the representations ${\cal R}_1,{\cal
R}_2$.

\begin{prop}
The representations $\{\Ps \colon \sigma\in S_n\}$ and $\{\Ms \colon \sigma \in S_n\}$ are
similar.
\end{prop}

\begin{proof}
The character of $\{\Ps \colon \sigma\in S_n\}$ is obviously given by $tr(\Ps)=\#\{i \colon
\sigma(i)=i\}$. We will show that $tr(\Ms)=\#\{i \colon \sigma(i)=i\}$ also, thereby
establishing that the representations are similar.

We first note that $tr({\cal M}_\sigma)=tr({\cal L}_\sigma)-h(\vw_{n+1})$, since every ``1'' in
$\vw_{n+1}$ is subtracted from exactly one diagonal position when we form $\Ms$ from $\Ls$ by
subtracting $\vw_{n+1}$ from each column. We show below that
 \begin{equation}\label{eq.traceL}
    tr(\Ls)=h(\vw_1)-1+\#\{i \colon \sigma(i)=i\},
 \end{equation}
so by Lemma~\ref{lem.2.4.2} we have
 $$tr(\Ms)=h(\vw_1)-1+\#\{i \colon \sigma(i)=i\}-h(\vw_{n+1})
                =\#\{i \colon \sigma(i)=i\},$$
which will conclude the proof.

Until this point we have found it convenient to think of $\Ls$ in terms of its columns. That is
not the natural viewpoint to take in proving Eq.~(\ref{eq.traceL}), however. The difference
between the $(j+1)$-st and $j$-th columns of $\Ls$ is
$\vd_{\sigma(j)}=\vec{e}_{i+1}-\vec{e}_i$; we think of this relationship as a ``1'' moving down
from the $i$-th row to the $(i+1)$-st row. In looking at the matrix $\Ls$ we see each ``1'' in
the first column continues across to the east, occasionally moving down a row (southeast), or
even `off' the bottom of the matrix. We call the path a ``1'' takes a {\em snake}.

\begin{xmp}
With $\sigma=[1,4,5,8,2,3,9,6,10,7]$, we find
 $$\Ls=\left(
 \begin{array}{ccccccccccc}
 {\bf1} & 0 & 0 & 0 & 0 & 0 & 0 & 0 & 0 & 0 & 0 \\
 0 & {\bf1} & {\bf1} & {\bf1} & {\bf1} & 0 & 0 & 0 & 0 & 0 & 0 \\
 0 & 0 & 0 & 0 & 0 & {\bf1} & 0 & 0 & 0 & 0 & 0 \\
 {\bf1} & {\bf1} & 0 & 0 & 0 & 0 & {\bf1} & {\bf1} & {\bf1} & {\bf1} & {\bf1} \\
 0 & 0 & {\bf1} & 0 & 0 & 0 & 0 & 0 & 0 & 0 & 0 \\
 0 & 0 & 0 & {\bf1} & {\bf1} & {\bf1} & {\bf1} & {\bf1} & 0 & 0 & 0 \\
 0 & 0 & 0 & 0 & 0 & 0 & 0 & 0 & {\bf1} & {\bf1} & 0 \\
 {\bf1} & {\bf1} & {\bf1} & {\bf1} & 0 & 0 & 0 & 0 & 0 & 0 & {\bf1} \\
 0 & 0 & 0 & 0 & {\bf1} & {\bf1} & {\bf1} & 0 & 0 & 0 & 0 \\
 0 & 0 & 0 & 0 & 0 & 0 & 0 & {\bf1} & {\bf1} & 0 & 0 \\
 \end{array} \right)$$
The matrix $\Ls$ has 3 snakes beginning in positions $(1,1)$, $(4,1)$, and $(8,1)$. The first
snake occupies the positions $(1,1)$, $(2,2)$, $(2,3)$, $(2,4)$, $(2,5)$, $(3,6)$, $(4,7)$,
$(4,8)$, $(4,9)$, $(4,10)$, and $(4,11)$.
\end{xmp}

Only the last snake moves off the bottom of $\Ls$; after all, $n$ only occurs once in a
permutation. The other snakes, of which there are $h(\vw_1)-1$, begin on or below the diagonal
and end on or above the diagonal. Thus each must intersect the diagonal at least once.
Moreover, each fixed point of $\sigma$ will keep a snake on a diagonal for an extra row. Thus,
$tr(\Ls)=h(\vw_1)-1+\#\{i \colon \sigma(i)=i\}$.
\end{proof}

We turn now to identifying the intertwining matrices, i.e., the matrices $\QQ$ such that
$\QQ^{-1} \Ms \QQ = \Ps$ for every $\sigma\in S_n$.

\begin{prop}\label{prop.similarityMatrices}
The $n\times n$ matrix $\QQ=(q_{ij})$ satisfies $\QQ^{-1} \Ms \QQ = \Ps$ for every $\sigma\in
S_n$ iff there are complex numbers $a,b$ with $(na+b)b^{n-1}\not=0$ and $q_{11}=a+b$,
$q_{1k}=a$, $q_{kk}=b$, $q_{k,k-1}=-b$ ($2\leq k\leq n$).
\end{prop}

\begin{proof}
We first note that it is sufficient to restrict $\sigma$ to a generating set of $S_n$. To see
this, let $S_n=\langle \sigma_1,\dots,\sigma_r\rangle$. If $\QQ$ satisfies $\QQ^{-1} \Ms \QQ =
\Ps$ for every $\sigma\in S_n$, then clearly $\QQ^{-1}{\cal M}_{\sigma_i} \QQ= {\cal
P}_{\sigma_i}$ ($1\leq i \leq r$). In the other direction, if $\QQ$ satisfies $\QQ^{-1}{\cal
M}_{\sigma_i} \QQ= {\cal P}_{\sigma_i}$ ($1\leq i \leq r$) and
$\sigma=\sigma_{i_1}\sigma_{i_2}\dots \sigma_{i_s}$, then
\begin{multline*}
    \QQ^{-1} \Ms \QQ =   \QQ^{-1}\left({\cal M}_{\sigma_{i_1}} {\cal M}_{\sigma_{i_2}}\dots
                        {\cal M}_{\sigma_{i_s}}\right) \QQ
                =   \prod_{j=1}^s \left( \QQ^{-1} {\cal M}_{\sigma_{i_j}} \QQ\right)\\
                =   \prod_{j=1}^s {\cal P}_{\sigma_{i_j}}
                =   {\cal P}_{\prod_{j=1}^s \sigma_{i_j}}
                =   {\cal P}_{\sigma}.
\end{multline*}

Thus we can restrict our attention to the transpositions $(k,k-1)$ ($2\leq k \leq n$). We
identified ${\cal M}_{(k,k-1)}$ in Lemma~\ref{lem.M_transposition} as ${\cal
M}_{(k,k-1)}=\II+\VV_k$, where $\VV_k$ is the matrix all of whose entries are zero, save the
$k$-th column, which is $\vec{e}_{k-1}-2\vec{e}_{k}+\vec{e}_{k+1}$.

We suppose that $\QQ=(q_{ij})$ satisfies ${\cal M}_{(k,k-1)}\QQ = \QQ {\cal P}_{(k,k-1)}$ to
find linear constraints on the unknowns $q_{ij}$. We will find that these constraints (for
$2\leq k \leq n$) are equivalent to $q_{11}=a+b$, $q_{1k}=a$, $q_{kk}=b$, $q_{k,k-1}=-b$
($2\leq k\leq n$). The determinant of $\QQ$ is easily seen to be $(na+b)b^{n-1}$, so that as
long as this is nonzero, $\QQ^{-1} \Ms \QQ = \Ps$ for every $\sigma \in S_n$.

We have ${\cal M}_{(k,k-1)}\QQ=\II \QQ+\VV_k \QQ$, so that ${\cal M}_{(k,k-1)}\QQ = \QQ
{\cal P}_{(k,k-1)}$ is equivalent to $\VV_k\QQ=\QQ({\cal P}_{(k,k-1)}-\II)$. This is a
convenient form since most entries in the matrices $\VV_k$ and ${\cal P}_{(k,k-1)}-\II$ are
zero. The entries of the product $\VV_k \QQ$ are 0 except for the $(k-1)$-st, $k$-th, and
$(k+1)$-st rows which are equal to the $k$-th row of $\QQ$, to $-2$ times the $k$-th row of
$\QQ$, and to the $k$-th row of $\QQ$, respectively. The entries of the product $\QQ({\cal
P}_{(k,k-1)}-\II)$ are 0 except for the $(k-1)$-st and $k$-th columns, which are equal to
the $k$-th column of $\QQ$ minus the $(k-1)$-st column of $\QQ$, and to the $(k-1)$-st
column of $\QQ$ minus the $k$-th column of $\QQ$, respectively. The entries which are zero
in one matrix or other lead to the families of equations $q_{kj}=0$ (for
$j\not\in\{k-1,k\}$) and $q_{jk}=q_{j,k-1}$ (for $j\not\in\{k-1,k,k+1\}$). The entries which
are non-zero in both products give the six equations
 $$\begin{pmatrix}
  q_{k,k-1}  &  q_{kk}  \\
 -2q_{k,k-1} & -2q_{kk} \\
  q_{k,k-1}  &  q_{kk}  \\
  \end{pmatrix}
  =
  \begin{pmatrix}
 q_{k-1,k}-q_{k-1,k-1} & q_{k-1,k-1}-q_{k-1,k} \\
   q_{kk}-q_{k,k-1}    &   q_{k,k-1}-q_{kk}    \\
 q_{k+1,k}-q_{k+1,k-1} & q_{k+1,k-1}-q_{k+1,k} \\
    \end{pmatrix},$$
which are equivalent to $q_{k-1,k-1}=q_{k-1,k}+q_{k,k}$, $q_{k,k-1}=-q_{kk}$, and
$q_{k+1,k-1}=q_{k+1,k}+q_{kk}$. Taking $q_{nn}=b$ and $q_{1n}=a$, the result follows.
\end{proof}

\begin{cor}\label{MsPs}
 $$\Ms =
 \begin{pmatrix}
 1  &         &   &        &         &   \\
 -1 &    1    &   &        & {\bf 0} &   \\
    &   -1    & 1 &        &         &   \\
    &         &   & \ddots &         &   \\
    & {\bf 0} &   &        & \ddots  &   \\
    &         &   &        &   -1    & 1 \\
 \end{pmatrix}
 \cdot \, \Ps \, \cdot
 \begin{pmatrix}
   1    &              &   &        &         &   \\
   1    &      1       &   &        & {\bf 0} &   \\
   1    &      1       & 1 &        &         &   \\
 \vdots &              &   & \ddots &         &   \\
 \vdots & {\bf 1}      &   &        & \ddots  &   \\
        & \hdots & \hdots  & \hdots &    1    & 1 \\
\end{pmatrix}.$$
\end{cor}

\begin{proof}
Set $a=0$, $b=1$ in Theorem~\ref{prop.similarityMatrices}. All that needs to be checked is that
$\QQ^{-1}$ is as claimed, i.e., the $n\times n$ matrix with ``1''s on and below and the
diagonal and ``0''s above the diagonal.
\end{proof}

We remark that, since it is easy to recover $\Ls$ from $\Ms$ (see the proof of
Lemma~\ref{lem.2.4.4}) this Corollary gives a simple method for computing the factors of length
$n$ of a Sturmian word with slope $\alpha$ given only the permutation ordering
$\fp{\alpha},\dots, \fp{n\alpha}$ (we don't even need $\alpha$). Also, we note that the matrix
with ``1''s on and below the diagonal is a summation operator, and its inverse is a difference
operator. If one could prove Corollary~\ref{MsPs} directly, this would provide a second proof
of Theorem~\ref{thm.Volume}.

\subsection{The Simplex $F_n(W)$}\label{sec.Sturmian.sub.Simplex}

Stolarsky \& Porta [personal communication] observed experimentally that $\Ma$ has determinant
$\pm 1$, and moreover that the roots of its characteristic polynomial are roots of unity. The
first observation was proved in the course of the proof of Theorem~\ref{thm.Volume} in
Subsection~\ref{sec.Sturmian.sub.Proof}. The second observation also follows from the fact that
$\Ma$ lies in a finite group.

We now summarize the results of this paper as they relate to $\Ma$ and $F_n(\alpha)$.

\begin{thm}\label{thm.summary}
Let $\alpha\in(0,1)$ be irrational, and $n\geq1$ an integer.
 \renewcommand{\theenumi}{\roman{enumi}}
 \begin{enumerate}
    \item
    The volume of the simplex $F_n(\alpha)$ is $\tfrac{1}{n!}$.
    \item
    $\det ( \Ma) = \pm 1$.
    \item
    $\Ma = {\cal M}_{\pan} = \QQ {\cal P}_{\pan} \QQ^{-1}$ (see
Theorem~\ref{prop.similarityMatrices} for the definition of $\QQ$);
    \item
    $\det({\cal M}_{2n}(\alpha)=\det({\cal M}_{2n+1}(\alpha))
        = \prod_{\ell=1}^n (-1)^{\floor{2\ell \alpha}}$.
    \item
    If $\pan(n)=n$, then
     $$\ord({\cal M}_{n-1}(\alpha))=\ord(\Ma)=\ord(\pan(1) \bmod{n}).$$
    \item
    If $\pan(1)=n$, then $$\ord({\cal M}_{n-1}(\alpha))=\ord(-\pan(n)
\bmod{n})$$ and $$\ord(\Ma)=\ord(-\pan(n) \bmod{gn}),$$ where $g$ is the smallest positive
integer such that $\gcd\left(n,\tfrac{\pan(n)+1}{g}\right)=1$.
  \end{enumerate}
\end{thm}

Items (i) and (ii) are proved in Subsection~\ref{sec.Sturmian.sub.Proof}. Item (iii) is a
combination of Theorem~\ref{prop.similarityMatrices} and
Proposition~\ref{prop.GoodDefinitions}. Items (iv), (v) and (vi) are immediate consequences of
Proposition~\ref{prop.GoodDefinitions}, the facts $\det(\Ms)=\sgn(\sigma)$ and
$\ord(\Ms)=\ord(\sigma)$, and Theorem~\ref{thm.PiAlphaN.properties}. They are included here for
the purpose of listing everything known about $\Ma$ in one place.

\section{Ordering Fractional Parts}\label{sec.FracParts}

\subsection{Statement of Results} Table~\ref{table.e} gives the sign and multiplicative order of
$\pi_{e,n}$ for $2\leq n \leq 136$. Visually inspecting the table, one quickly notices that
$\sgn(\pi_{e,2n})=\sgn(\pi_{e,2n+1})$ for all $n$, and that $\ord(\pi_{e,n})$ is surprisingly
small for $n=70,71,109,110$. The first several convergents to $e$ are $2, 3, \frac{8}{3},
\frac{11}{4}, \frac{19}{7}, \frac{87}{32}, \frac{106}{39}, \frac{193}{71}$; the value 71 is the
denominator of a convergent, and $110=71+39$ is the sum of two denominators. Note also that for
some values of $n$, $\ord(\pi_{e,n})$ is extraordinarily large, e.g.,
$\ord(\pi_{e,123})=22383900$. None of these observations are peculiar to the irrational $e$.
Some of these observations are explained by Theorem~\ref{thm.PiAlphaN.properties} below, and
the others remain conjectural.

\begin{table}[b!]
\caption{Algebraic properties of $\pan$ with $\alpha=e$ and $2\leq n \leq 136$}\label{table.e}
\vspace{3mm}
    {\tiny
\begin{center}
\begin{tabular}{|rcc|rcc|rcc|}
\hline
 $n$ & $\sgn(\pi_{\alpha,n})$ & $\ord(\pi_{\alpha,n})$ & $n$ & $\sgn(\pi_{\alpha,n})$ & $\ord(\pi_{\alpha,n})$ & $n$ & $\sgn(\pi_{\alpha,n})$ & $\ord(\pi_{\alpha,n})$ \\
\hline
     &                        &                        &     &                        &                        &     &                        &                        \\
   2 &           -1           &           2            &  47 &           -1           &           44           &  92 &           1            &          2107          \\
   3 &           -1           &           2            &  48 &           -1           &          540           &  93 &           1            &         13244          \\
   4 &           -1           &           2            &  49 &           -1           &          120           &  94 &           -1           &         18810          \\
   5 &           -1           &           4            &  50 &           1            &          680           &  95 &           -1           &         20034          \\
   6 &           -1           &           6            &  51 &           1            &          1848          &  96 &           -1           &          3348          \\
   7 &           -1           &           6            &  52 &           -1           &           50           &  97 &           -1           &         11256          \\
   8 &           1            &           7            &  53 &           -1           &           90           &  98 &           -1           &          1702          \\
   9 &           1            &           6            &  54 &           -1           &          962           &  99 &           -1           &          188           \\
  10 &           -1           &           10           &  55 &           -1           &          1848          & 100 &           1            &          957           \\
  11 &           -1           &           10           &  56 &           -1           &          588           & 101 &           1            &          2100          \\
  12 &           -1           &           12           &  57 &           -1           &          276           & 102 &           -1           &          102           \\
  13 &           -1           &           36           &  58 &           1            &          165           & 103 &           -1           &          2052          \\
  14 &           -1           &           40           &  59 &           1            &          1260          & 104 &           -1           &          1950          \\
  15 &           -1           &           14           &  60 &           -1           &          1848          & 105 &           -1           &          1260          \\
  16 &           1            &           15           &  61 &           -1           &          2040          & 106 &           -1           &          5964          \\
  17 &           1            &           3            &  62 &           -1           &           62           & 107 &           -1           &         13860          \\
  18 &           1            &           3            &  63 &           -1           &         15640          & 108 &           1            &         54366          \\
  19 &           1            &           15           &  64 &           1            &          2040          & 109 &           1            &           10           \\
  20 &           1            &           77           &  65 &           1            &          424           & 110 &           -1           &           10           \\
  21 &           1            &           12           &  66 &           -1           &          966           & 111 &           -1           &          2310          \\
  22 &           -1           &           12           &  67 &           -1           &          1476          & 112 &           -1           &          720           \\
  23 &           -1           &           12           &  68 &           -1           &           56           & 113 &           -1           &          3738          \\
  24 &           1            &           4            &  69 &           -1           &          232           & 114 &           1            &          1938          \\
  25 &           1            &           4            &  70 &           -1           &           14           & 115 &           1            &          112           \\
  26 &           1            &           36           &  71 &           -1           &           14           & 116 &           -1           &         92820          \\
  27 &           1            &           48           &  72 &           1            &          6840          & 117 &           -1           &         11220          \\
  28 &           1            &           6            &  73 &           1            &          406           & 118 &           -1           &          5520          \\
  29 &           1            &           24           &  74 &           -1           &          390           & 119 &           -1           &         60060          \\
  30 &           -1           &          210           &  75 &           -1           &          780           & 120 &           -1           &         14280          \\
  31 &           -1           &           4            &  76 &           -1           &          192           & 121 &           -1           &          1680          \\
  32 &           -1           &           4            &  77 &           -1           &          228           & 122 &           1            &          6240          \\
  33 &           -1           &          180           &  78 &           -1           &         32130          & 123 &           1            &        22383900        \\
  34 &           -1           &          210           &  79 &           -1           &          390           & 124 &           -1           &         820820         \\
  35 &           -1           &          420           &  80 &           1            &          630           & 125 &           -1           &         215460         \\
  36 &           1            &          120           &  81 &           1            &           72           & 126 &           -1           &          9360          \\
  37 &           1            &           37           &  82 &           1            &          2728          & 127 &           -1           &         17160          \\
  38 &           -1           &           12           &  83 &           1            &          6138          & 128 &           1            &         68640          \\
  39 &           -1           &           12           &  84 &           1            &          152           & 129 &           1            &         47888          \\
  40 &           -1           &           40           &  85 &           1            &          6669          & 130 &           -1           &          7276          \\
  41 &           -1           &          1980          &  86 &           -1           &         31920          & 131 &           -1           &          508           \\
  42 &           -1           &          414           &  87 &           -1           &          400           & 132 &           -1           &          6720          \\
  43 &           -1           &           42           &  88 &           1            &          192           & 133 &           -1           &          4914          \\
  44 &           1            &          580           &  89 &           1            &         14616          & 134 &           -1           &          1560          \\
  45 &           1            &          168           &  90 &           1            &         18585          & 135 &           -1           &         11752          \\
  46 &           -1           &          1120          &  91 &           1            &         25080          & 136 &           1            &          3045          \\
\hline
\end{tabular}
\end{center}
    }
\end{table}

As in Section~\ref{sec.Sturmian}, we make frequent use of Knuth's notation:
 $$\tf{Q}=\begin{cases} 1 & \text{$Q$ is true;} \\ 0 & \text{$Q$ is false.}\end{cases}$$

Lemma~\ref{lem.Sos}, giving $\pan$ in terms of only $\pan(1)$, $\pan(n)$, and $n$, is proved
by V. T. S{\'o}s in~\cite{Sos1957}. Her method of proof is similar to our proof of
Lemma~\ref{Brepresentation} below. The lemma is also proved---in terse English---in
\cite{MR36:114}. We will derive Theorem~\ref{thm.PiAlphaN.properties} from S{\'o}s's Lemma.
\begin{lem}[S{\'o}s]\label{lem.Sos}
Let $\al$ be irrational, $n$ a positive integer, and $\pi=\pan$. Then
  $$\pi(k+1)=\pi(k)+\pi(1)\,\tf{\pi(k)\le \pi(n)}-\pi(n)\,\tf{n<\pi(1)+\pi(k)}$$
for $1\leq k < n$.
\end{lem}

The surprising Three-Distance Theorem is an easy corollary: If $\alpha$ is irrational, the
$n+2$ points $0,\fp{\alpha},\fp{2\alpha},\dots,\fp{n\alpha},1$ divide the interval $[0,1]$ into
$n+1$ subintervals which have at most 3 distinct lengths. Alessandri \&
Berth{\'e}~\cite{AlessandriBerthe1998} give an excellent and up-to-date survey of generalizations
of the Three-Distance Theorem.

The primary goal of this section is to prove Theorem~\ref{thm.PiAlphaN.properties}, which
refines Theorem~\ref{thm.Order}, and to prove Theorem~\ref{thm.Sign}.
Corollary~\ref{phiorder} is of independent interest.

\begin{thm}
\label{thm.PiAlphaN.properties} Let $\alpha \not\in \Q$ and $n\in\Z^+$.
    \renewcommand{\theenumi}{\roman{enumi}}
    \begin{enumerate}
        \item If $\pan(n)=n$, then $\ord(\pi_{\alpha,n-1}) = \ord(\pan) =
                                                \ord(\pan(1) \bmod{n})$.
        \item If $\pan(1)=n$, then $\ord(\pi_{\alpha,n-1})=\ord(-\pan(n)
                    \bmod{n})$, and $$\ord(\pan)=\ord(-\pan(n)\bmod{gn}),$$ where
                    $g$ is the least positive integer such that
                    $\gcd\left(n,\frac{\pan(n)+1}{g}\right)=1$.
    \end{enumerate}
\end{thm}

\subsection{The Multiplicative Order of the Permutation}

Theorem~\ref{thm.PiAlphaN.properties} follows from S{\'o}s's Lemma.

\begin{proofof}{Theorem~\ref{thm.PiAlphaN.properties}(i)}
Suppose that $\pan(n)=n$. We have
 $$\pan(k)=\begin{cases} \pi_{\alpha,n-1}(k) & 1\leq k \leq n-1; \\ k & k=n,
    \end{cases}$$
and so obviously $\ord(\pi_{\alpha,n-1})=\ord(\pan)$. We show that the length of every orbit
divides $\ord(\pan(1) \bmod{n})$, and that the length of the orbit of 1 is exactly
$\ord(\pan(1) \bmod{n})$. From S{\'o}s's Lemma (Lemma~\ref{lem.Sos}), we have in this case for
$1\leq k <n$ the congruence $\pan(k+1) \equiv \pan(k)+\pan(1) \pmod{n}.$ By induction, this
gives $\pan(k)\equiv k \pan(1) \pmod{n}$ for $1\leq k \leq n$. Thus, the orbit of the point $k$
is
 $$k,\, k\pan(1),\, k\pan(1)^2,\, k\pan(1)^3,\, \dots$$
The length of the orbit of $k$ divides $\ord(\pan(1) \bmod{n})$, and in particular the orbit of
1 has length equal to $\ord(\pan(1) \bmod{n})$.
\end{proofof}

\begin{proofof}{Theorem~\ref{thm.PiAlphaN.properties}(ii)}
Suppose that $\pan(1)=n$. S{\'o}s's Lemma gives $$\pan(k)\equiv (1-k)\pan(n) \pmod{n}.$$ For $1\leq
k <n$ we have $$\pi_{\alpha,n-1}(k)=\pan(k+1)\equiv(1-(k+1))\pan(n)=-k\pan(n) \pmod{n}.$$ Thus
for $r\geq 1$ we have $\pi_{\alpha,n-1}^r(k)\equiv k (-\pan(n))^r \pmod{n}$. The length of the
orbit of $k$ under $\pi_{\alpha,n-1}$ divides $\ord(-\pan(n) \bmod{n})$, and in particular the
orbit of 1 has length equal to $\ord(-\pan(n) \bmod{n})$. This proves that
$\ord(\pi_{\alpha,n-1})=\ord(-\pan(n) \bmod{n})$.

For notational convenience, set $x=-\pan(n)$. Now from
 $$\pan(k)\equiv (1-k)\pan(n)=(k-1)x \pmod{n}$$
it is readily seen by induction that for $R\geq 0$ we have
 \begin{equation}\label{Rimage}
   \pan^{R}(k) \equiv k x^R-
                \left(x+x^2+\dots+x^R\right) \pmod{n}.
 \end{equation}
Let $r$ be the least positive integer for which $\forall k \left(\pan^r(k)=k\right)$, i.e., let
$r$ be the least common multiple of the length of the orbits of $\pan$. We must show that
$r=\ord(x\bmod{gn})$.

Define the integer $\gamma$ by $g\gamma=x-1$, and note that $\gcd(\gamma,n)=1$. We may
rearrange Eq.~(\ref{Rimage}), setting $R=r$, as
 \begin{equation}\label{rproperty}
    \frac{x^r-1}{x-1} \equiv (k-1)(x^r-1) \pmod{n},
 \end{equation}
the division being real, not modular. With $k=1$, Eq.~(\ref{rproperty}) becomes $0\equiv
\tfrac{x^r-1}{x-1}=\tfrac{x^r-1}{g\gamma} \pmod{n}$, which holds iff $\tfrac{x^r-1}{g} \equiv 0
\pmod{n}$. This, in turn, holds iff there is an integer $\beta$ with $\beta n
=\frac{x^r-1}{g}$, i.e., $\beta n g = x^r-1$. Thus $r$ is a multiple of $\ord( x \bmod{gn})$,
and in particular $r\geq \ord( x \bmod{gn})$.

We claim that $\forall k \left( \pan^{\ord( x \bmod{gn})}(k)=k\right)$, so that $r \leq \ord( x
\bmod{gn})$, which will conclude the proof. We have
 $$
    \frac{x^{\ord( x \bmod{gn})}-1}{x-1}
        =\frac{\beta g n}{g\gamma}
        =\frac{\beta n}{\gamma}
 $$
and
 $$
    (k-1)(x^{\ord( x \bmod{gn})}-1) \equiv 0 \pmod{n},$$
so that substituting $r=\ord( x \bmod{gn})$ into Eq.~(\ref{rproperty}) we write $\frac{\beta
n}{\gamma} \equiv 0 \pmod{n}$, which holds for all $k$ since $\gcd(\gamma,n)=1$. Thus $r \leq
\ord( x \bmod{gn})$.
\end{proofof}

For quadratic irrationals one can easily identify the convergents and intermediate fractions
and, if the denominators have enough structure, explicitly compute $\ord(\pan)$ when $n$ or
$n+1$ is such a denominator. We have for example

\begin{cor}\label{phiorder}
Let $\phi=\tfrac{\sqrt{5}-1}{2}$, and $f_n$ be the $n$-th Fibonacci number. Then for $n\geq2$,
$$\ord(\pi_{\phi,-1+f_{2n}})=\ord(\pi_{\phi,f_{2n}})=2$$ and
$$\ord(\pi_{\phi,-1+f_{2n+1}})=\ord(\pi_{\phi,f_{2n+1}})=4.$$
\end{cor}

\begin{proof}
From the continued fraction expansion of $\phi$ we know that $\pi_{\phi,f_{2n}}(f_{2n})=f_{2n}$
and $\pi_{\phi,f_{2n}}(1)=f_{2n-1}$. The identity $f_{2n-1}^2-f_{2n}f_{2n-2}=1$ shows that
$\ord(f_{2n-1} \bmod{f_{2n}})=2$ for $n\geq2$, whence by
Theorem~\ref{thm.PiAlphaN.properties}(i) we have
$\ord(\pi_{\phi,-1+f_{2n}})=\ord(\pi_{\phi,f_{2n}})=2$.

The continued fraction expansion of $\phi$ also tells us $\pi_{\phi,f_{2n+1}}(1)=f_{2n+1}$ and
\linebreak $\pi_{\phi,f_{2n+1}}(f_{2n+1})=f_{2n}$. The identity $f_{2n}^2-f_{2n+1}f_{2n-1}=-1$
shows $\ord(-f_{2n} \bmod{f_{2n+1}})=4$ for $n\geq2$, whence by
Theorem~\ref{thm.PiAlphaN.properties}(ii), $\ord(\pi_{\phi,-1+f_{2n+1}})=4$.

We defined $g$ to be the least positive integer such that
$\gcd\left(f_{2n+1},\frac{f_{2n}+1}{g}\right)=1$. In particular, $g|f_{2n}+1$, and so the
identity
 $$ (-f_{2n})^4-1 = f_{2n-1}(f_{2n}-1)f_{2n+1}(f_{2n}+1)
        \equiv 0 \pmod{g f_{2n+1}}$$
shows that $\ord(-f_{2n} \bmod{g f_{2n+1}})$ divides 4, and the identity
$(-f_{2n})^2-f_{2n+1}f_{2n-1}=-1$ shows that $$\ord(-f_{2n} \bmod{g f_{2n+1}}) \geq
\ord(-f_{2n} \bmod{f_{2n+1}}) =4,$$ whence $\ord(-f_{2n} \bmod{g f_{2n+1}})=4$ and by
Theorem~\ref{thm.PiAlphaN.properties}(ii), $\ord(\pi_{\phi,f_{2n+1}})=4$.
\end{proof}

\subsection{The Sign of the Permutation}

Recall that a $k$-cycle is even exactly if $k$ is odd. We define $\rho(n,k)$ to be the
$(n-k)$-cycle
    $$\rho(n,k):= (n,n-1,\dots,k+1) = (n,n-1)(n-1,n-2)\cdots (k+2,k+1).$$

Define also $$B_\alpha(k):=\#\{q \colon 1\leq q < k,\, \fp{q\alpha}<\fp{k\alpha}\},$$ which
counts the integers in $[1,k)$ that are `better' denominators for approximating $\alpha$
from below. Clearly,
 \begin{align*}
 \pan       &=   \pi_{\alpha,n-1}\,\rho(n,\B(n)) \\
            &=   \rho(1,\B(1))\,\rho(2,\B(2))\,\dots\,\rho(n,\B(n)) \\
            &=   \prod_{k=1}^n \rho(k,\B(k))
 \end{align*}
so that $\pan$ is the product of $\sum_{k=1}^n (k-\B(k)-1)$ transpositions. We will show that
for $k$ odd, $\B(k)\equiv 0 \pmod{2}$, which will be used to demonstrate that
$\sgn(\pi_{\alpha,2n})=\sgn(\pi_{\alpha,2n+1})$.

Our proof of Lemma~\ref{Brepresentation} is similar in spirit to S{\'o}s's proof of
Lemma~\ref{lem.Sos}.

\begin{lem}\label{Brepresentation}
For $k\geq 3$ and $0<\alpha<1/2$, $\alpha$ irrational, $\B(k)+B_{1-\alpha}(k)=k-1$, and
 $$\B(k) -2\B(k-1)+\B(k-2) =
   \begin{cases}
    1-k,     & \fp{k\alpha}  \in [0,\alpha); \\
    k-1,     & \fp{k\alpha}  \in [\alpha,2\alpha); \\
    0,       & \fp{k\alpha}  \in [2\alpha,1).
  \end{cases}
  $$
\end{lem}

\begin{proof}
Observe that $0<\fp{q\alpha} < \fp{k\alpha}$ iff $\fp{k(1-\alpha)}<\fp{q(1-\alpha)}<1$, so that
$q$ with $1\leq q<k$ is in either the set $\{q: 1\leq q<k, \fp{q\alpha}<\fp{k\alpha}\}$ or in
the set $\{ q: 1\leq q<k, \fp{q(1-\alpha)}<\fp{k(1-\alpha)}\}$, and is not in both. Thus,
$\B(k)+B_{1-\alpha}(k)=k-1$.

We think of the points $0,\fp{\alpha},\dots,\fp{k \alpha}$ as lying on a circle with
circumference 1, and labeled $P_0, P_1, \dots, P_k$, respectively, i.e., $P_j :=
\tfrac{1}{2\pi} e^{2\pi j \alpha \sqrt{-1}}=\frac{1}{2\pi}e^{2\pi \fp{j \alpha} \sqrt{-1}}$.
``The arc $\overline{P_iP_j}$'' refers to the half-open counterclockwise arc from $P_i$ to
$P_j$, containing $P_i$ but not $P_j$. We say that three distinct points $A,B,C$ are {\em in
order} if $B \not\in \overline{CA}$. We say that $A,B,C,D$ are in order if both $A,B,C$ and
$C,D,A$ are in order. Essentially, if when moving counter-clockwise around the circle starting
from $A$, we encounter first the point $B$, then $C$, then $D$, and finally $A$ (again), then
$A,B,C,D$ are in order.

By rotating the circle so that $P_i \mapsto P_{i+1}$ ($0\leq i \leq k$), we find that each $P$
on the arc $\overline{P_{k-2}P_{k-1}}$ is rotated onto a $P$ on the arc
$\overline{P_{k-1}P_k}$. Specifically, the number of $P_0,P_1,\dots,P_{k-2}$ on
$\overline{P_{k-2}P_{k-1}}$ is the same as the number of $P_1,P_2,\dots,P_{k-1}$ on
$\overline{P_{k-1}P_k}$. Set
 $$X:=\{P_0,P_1,\dots,P_{k-2}\} \quad\text{and}\quad
    Y:=\{P_1,P_2,\dots,P_{k-1}\},$$
so that what we have observed is
 \begin{equation}\label{dagger}
    \left| X \cap \overline{P_{k-2}P_{k-1}} \right| =
    \left| Y \cap \overline{P_{k-1}P_k} \right|.
 \end{equation}
Also, we will use
 $$ \B(k) = \left| Y \cap \overline{P_0P_k}\right|.$$

Now, first, suppose that $\fp{k\alpha} \in [0,\alpha)$, so that the points $P_0, P_k, P_{k-2},
P_{k-1}$ are in order on the circle. We have
  \begin{align*}
    X \cap \overline{P_{k-2}P_{k-1}}
        &=  X \cap \left( \overline{P_0P_{k-1}} \setminus
                \overline{P_0P_{k-2}}\right)\\
        &=  \left( X \cap \overline{P_0P_{k-1}}\right) \setminus
                \left(X \cap \overline{P_0P_{k-2}}\right) \\
    \left|X \cap \overline{P_{k-2}P_{k-1}}\right|
        &=  \left|\left( X \cap \overline{P_0P_{k-1}}\right) \right|\,-\,
                \left| \left(X \cap \overline{P_0P_{k-2}}\right)\right| \\
        &=  \left(\B(k-1)+1\right)-\left(\B(k-2)+1\right)\\
        &=  \B(k-1)-\B(k-2),
  \end{align*}
and similarly
  \begin{align*}
    Y \cap \overline{P_{k-1}P_k}
        &=  \left( Y \cap \overline{P_{k-1}P_0}\right) \cup
                \left( Y \cap \overline{P_0P_k}\right)\\
        &=  \left( Y \setminus \left( Y \cap \overline{P_0P_{k-1}} \right) \right)
                \cup \left( Y \cap \overline{P_0P_k} \right) \\
    \left| Y \cap \overline{P_{k-1}P_k} \right|
        &=  (|Y|-\left| Y \cap \overline{P_0P_{k-1}} \right|)\,+\,
                \left|Y \cap \overline{P_0P_k} \right| \\
        &=  (k-1-\B(k-1))+\B(k)\\
        &= \B(k)-\B(k-1)-(1-k),
  \end{align*}
so that Eq.~(\ref{dagger}) becomes $\B(k-1)-\B(k-2)=\B(k)-\B(k-1)-(1-k)$, as claimed in the
statement of the theorem.

Now suppose that $\fp{k\alpha} \in [\alpha,2\alpha)$, so that the points $P_0, P_{k-1}, P_k,
P_{k-2}$ are in order. By arguing as in the above case, we find
$$X \cap \overline{P_{k-2}P_{k-1}}=\left(X\setminus\left(X\cap
\overline{P_0P_{k-2}}\right)\right)\cup \left(X\cap \overline{P_0P_{k-1}}\right),$$ and so
$\left|X\cap\overline{P_{k-2}P_{k-1}}\right|=\B(k-1)-\B(k-2)+(k-1)$. Likewise,
$$Y \cap \overline{P_{k-1}P_k}=\left(Y\cap\overline{P_0P_k}\right)\setminus
\left(Y\cap\overline{P_0P_{k-1}}\right)$$ so that $\left|Y \cap \overline{P_{k-1}P_k}
\right|=\B(k)-\B(k-1)$. Thus, in this case Eq.~(\ref{dagger}) becomes
$\B(k-1)-\B(k-2)+(k-1)=\B(k)-\B(k-1)$, as claimed in the statement of the theorem.

Finally, suppose that $\fp{k\alpha} \in [2\alpha,1)$, so that the points
$P_0,P_{k-2},P_{k-1},P_k$ are in order. We find
 $$X\cap\overline{P_{k-2}P_{k-1}}=\left( X \cap \overline{P_0P_{k-1}}\right)
 \setminus \left(X \cap \overline{P_0P_{k-2}}\right)$$
and $\left|X\cap\overline{P_{k-2}P_{k-1}}\right|=\B(k-1)-\B(k-2)$. Also,
 $$Y\cap\overline{P_{k-1}P_k} = \left(Y\cap\overline{P_0P_k}\right)\setminus
 \left(Y\cap\overline{P_0P_{k-1}}\right)$$
and so $\left|Y\cap\overline{P_{k-1}P_k}\right|=\B(k)-\B(k-1)$. As claimed, we have
$\B(k-1)-\B(k-2)=\B(k)-\B(k-1)$.
\end{proof}

Lemma~\ref{SignCorollary} makes explicit the connection between Lemma~\ref{Brepresentation},
arithmetic properties of $B_\alpha$, and the permutation $\pan$.

\begin{lem} \label{SignCorollary}
Let $\alpha\in(0,1)$ be irrational and $n\in\Z^+$. If $k$ is odd, then $\B(k)$ is even. If
$k$ is even, then $\B(k) \equiv \floor{k\alpha}+1 \pmod{2}$.
\end{lem}

\begin{proof}
By Lemma~\ref{Brepresentation}, $\B(k)+B_{1-\alpha}(k)=k-1$. Thus for odd $k$,
$\B(k)+B_{1-\alpha}(k)$ is even, and so either both $\B(k)$ and $B_{1-\alpha}(k)$ are even
or both are odd. This means that, for odd $k$, without loss of generality we may assume that
$0<\alpha<\tfrac12$.

Reducing the recurrence relation in Lemma~\ref{Brepresentation} modulo 2, under the
hypothesis that $k$ is odd, we find that $\B(k) \equiv \B(k-2) \pmod{2}$. Since $\B(1)=0$,
we see that $\B(k)\equiv0 \pmod{2}$ for all odd $k$.

Now suppose that $k$ is even and $0<\alpha<\tfrac12$. The recurrence relation in
Lemma~\ref{Brepresentation} reduces to
 $$\B(k)+\B(k-2) \equiv \tf{\fp{k\alpha} \in [0,2\alpha)} \pmod{2}.$$
Set $\beta=2\alpha$, $k=2\ell$ and $\BB(i)=\B(2i)$. We have
 \begin{align*}
    \B(k)=\B(2\ell)
                &=      \BB(\ell) \\
                &\equiv \BB(\ell-1)+\tf{\fp{\ell\beta} \in [0,\beta)}
                                \pmod{2}\\
                &\equiv \BB(1)+\sum_{i=2}^{\ell}
                         \tf{\fp{i\beta} \in [0,\beta)} \pmod{2}\\
                &=      \B(2)+\floor{\ell\beta} \\
                &=      1+\floor{2\ell\alpha} = 1+\floor{k\alpha},
 \end{align*}
since $\sum_{i=2}^{\ell} \tf{\fp{i\beta} \in [0,\beta)}$ (with $\beta\in(0,1)$) counts the
integers in the interval $(\beta,\ell\beta]$. This proves the lemma for $0<\alpha<\tfrac12$
and $k$ even.

Now suppose that $k$ is even and $\tfrac12<\alpha<1$. By Lemma~\ref{Brepresentation}, $\B(k)
= k-1-B_{1-\alpha}(k) \equiv 1+B_{1-\alpha}(k) \pmod{2}$. By the argument (for
$0<\alpha<\tfrac12$) given above, $B_{1-\alpha}(k) \equiv 1+\floor{k(1-\alpha)}\pmod{2}$. We
have
 \begin{align*}
  \B(k) &\equiv \floor{k(1-\alpha)}\pmod{2} \\
        &=      k-\floor{k\alpha}-\fp{k(1-\alpha)}-\fp{k\alpha} \\
        &\equiv \floor{k\alpha}+1 \pmod{2},
 \end{align*}
where we have used $\fp{k(1-\alpha)}+\fp{k\alpha}=1$ (since $\alpha$ is irrational).
\end{proof}

\begin{proofof}{Theorem~\ref{thm.Sign}}
We have $\pi_{\alpha,2n+1}=\pi_{\alpha,2n} \, \rho(2n+1,\B(2n+1))$. The permutation
$\rho(2n+1,\B(2n+1))$ is the product of $2n+1-\B(2n+1)-1$ transpositions, which is an even
number by Lemma~\ref{SignCorollary}. Thus $\sgn(\pi_{\alpha,2n+1})=\sgn(\pi_{\alpha,2n})$.

By Lemma~\ref{SignCorollary}, for all integers $k$
 $$k-\B(k)+1 \equiv
 \left\{%
    \begin{array}{ll}
        0, & \hbox{$k$ odd;} \\
        \floor{k\alpha}, & \hbox{$k$ even.} \\
    \end{array}%
 \right.
 $$
Since the sign of the permutation $\rho(k,\B(k))$ is $(-1)^{k-\B(k)}$ and $(-1)^{2n}=1$, we
have
 \begin{multline*}
 \sgn(\pi_{\alpha,2n}) =  (-1)^{2n} \sgn\left( \prod_{k=1}^{2n}
                                                     \rho(k,\B(k))\right)
            = \prod_{k=1}^{2n} (-1)\sgn\left( \rho(k,\B(k)) \right) \\
            = \prod_{k=1}^{2n} (-1)(-1)^{k-\B(k)}
            = \prod_{k=1}^{2n} (-1)^{k-\B(k)+1}
            = \prod_{\ell=1}^n (-1)^{\floor{2\ell\alpha}}.
 \end{multline*}
\end{proofof}

\section{Unanswered Questions}\label{sec.Questions}

The most significant question we have been unable to answer is {\em why} the matrices formed
from Sturmian words in Section 2 lie in a common representation of $S_n$. We made two choices
with little motivation: we ordered the factors anti-lexicographically; and we subtracted the
last factor from the others. What happens if we order the factors differently, or subtract the
second factor from the others? Understanding why the structure revealed in Section 2 exists
might allow us to predict other phenomena.

Lucas Wiman [personal communication] has proved that $$\left\{{\cal M}_{n-1}\left(\tfrac
cn\right) \colon \gcd(c,n)=1\right\}$$ is isomorphic to the multiplicative group modulo $n$,
and asks if this is the largest subset of $\{{\cal M}_{n-1}(\alpha)\colon 0<\alpha<1\}$ that
is a group.

We have shown that the volume of the simplex $F_n(\alpha)$ is independent of $\alpha$. When are
two such simplices actually congruent?

\begin{cnj}
For $\alpha,\beta\in(0,1)$ be irrational, $F_n(\alpha) \cong F_n(\beta)$ (as simplices) iff
$F_n(\alpha)=F_n(\beta)$ or $F_n(\alpha)=F_n(1-\beta)$.
\end{cnj}

We have verified this conjecture for $n\leq 20$ by direct computation. At least, can any
such simplex be cut and reassembled (in the sense of Hilbert's 3rd Problem: see
Eves~\cite{Eves} for the basic theory and Sydler~\cite{Sydler} for a complete
characterization) into the shape of another?

While we have identified the matrix $\Ma$, there are many interesting questions that we remain
unable to answer. We don't believe that there is a bound on $\sum_{n=1}^N
\det(\Ma)=\sum_{n=1}^N \sgn(\pan)$ that is independent of $N$, but this sum must grow very
slowly. We suspect that
 $$\left| \sum_{n=1}^N \sgn(\pan) \right| \ll \log N$$
for almost all $\alpha$. For example, with $\alpha=\tfrac{\sqrt{5}-1}{2}$ and $N< e^{13}\approx
442413$,
 $$\left| \sum_{n=1}^N \sgn(\pan) \right| < 10.$$

In Section 3.1 we showed that $\ord(\pan)$ is regularly extremely small relative to the average
order of a permutation on $n$ symbols. For each irrational $\alpha$, are there infinitely many
values of $n$ for which $\ord(\pan)$ is exceptionally large?

One might hope for an explicit formula for $\ord(\pan)$ in terms of the base-$\alpha$ Ostrowski
expansion of $n$, but this seems to be extremely difficult. Are metric results more
approachable? Specifically, what is the distribution of $\ord(\pan)$ and $\sgn(\pan)$ for
$\alpha$ taken uniformly from $(0,1)$?

\begin{figure}
\begin{center}
    \begin{picture}(320,200)
        \put(330,12){$n$}
        \put(16,210){$I(n)$}
        \includegraphics{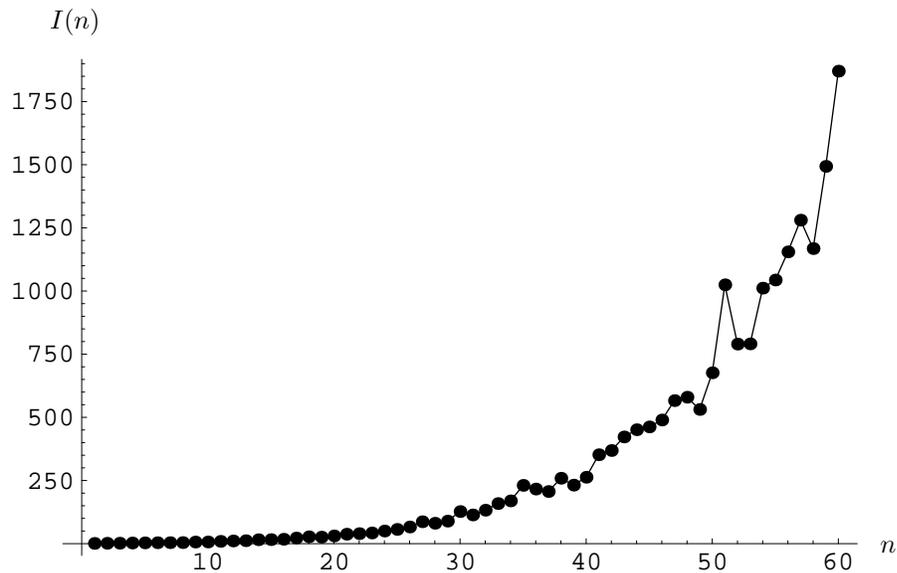}
    \end{picture}
\caption{$I(n)$ for $1\le n \le 60$} \label{pic.I(n)}
\end{center}
\end{figure}
Since $\ord(\pan)$ appears to vary wildly, it may be advantageous to consider its average
behavior. Can one give an asymptotic expansion of $\sum_{n=1}^N \ord(\pan)$? What can be
said about $I(n):=\int_0^1 \ord(\pan)\,d\alpha$? Surprisingly, although $I(n)$ seems to be
rapidly increasing, it is not monotonic, e.g., $I(35)>I(36)>I(37)$. Is this the `law of
small numbers,' or are there infinitely many values of $n$ for which $I(n+1)>I(n)$? We are
not aware of any non-trivial bounds, upper or lower, on $I(n)$. Figure~\ref{pic.I(n)} shows
$I(n)$ for $1\leq n \leq 60$.

$\B$ is an interesting function in its own right. We gave a formula for $\B(n) \pmod{2}$ in
Lemma~\ref{SignCorollary}. Is it possible to give a nice formula for $\B(n)$ for other
moduli? It seems likely that there are arbitrarily large integers which are not in the range
of $\B$, although we have been unable to show that it does not contain {\em all} nonnegative
integers. For example, $B_{(\sqrt{5}-3)/2}(k) \not= 3, B_{\sqrt{2}}(k)\not= 7$ and
$B_{e^{-1}}(k)\not=23$ for $k\leq 10^7$. It is perhaps noteworthy that the least $k$ for
which $B_{e^{-1}}(k)=25$ is $k=22154$, a reminder that $\B$ can take new, small values even
at large $k$. From the theory of continued fractions we know that there are infinitely many
$k$ for which $\B(k)=0$. Is there an $x\not=0$ and irrational $\alpha$ for which there are
infinitely many $k$ such that $\B(k)=x$?

\bigskip{\Large \noindent \bf Acknowledgements}\medskip

\noindent I am grateful to Horacio Porta \& Kenneth B. Stolarsky for bringing $\Ma$ to my
attention, and for making their conjectures and data available. I also wish to thank my
thesis advisor, Kenneth B. Stolarsky, for his guidance and careful reading, and Lucas Wiman
for his many probing questions and suggestions.


\end{document}